\documentclass[11pt]{article}
\usepackage[letterpaper,hmargin=1.06in,top=.96in,bottom=.98in]{geometry}
\usepackage[T1]{fontenc}
\usepackage{lmodern,amsmath,amssymb,amsthm,microtype}
\usepackage{tikz}
\usetikzlibrary{arrows.meta}
\usepackage[hidelinks]{hyperref}
\hypersetup{pdftitle={Large blow-ups of uniform hypergraphs},pdfauthor={Dhruv Mubayi}}
\newtheorem{theorem}{Theorem}
\newtheorem{corollary}[theorem]{Corollary}
\newtheorem{lemma}[theorem]{Lemma}
\newcommand{\cK}{K_h^{(k)}}
\newcommand{\cX}{\mathcal X}
\newcommand{\aut}{\operatorname{Aut}}
\newcommand{\eps}{\varepsilon}
\newcommand{\F}{F_2^{(k)}}
\title{Large blow-ups of uniform hypergraphs}
\author{Dhruv Mubayi\thanks{Department of Mathematics, Statistics and
		Computer Science, University of Illinois, Chicago, IL 60607.
		Email: mubayi@uic.edu. Research partially supported by NSF Awards DMS-2153576 and DMS-2552740.}}
\date{November 23, 2026}

\begin{document}
	\maketitle
	
	\begin{abstract}
		
		Let \(k\ge2\), let \(H\) be a fixed \(k\)-uniform hypergraph on \(h\) vertices
		with at least one edge, and let \(a>0\).
		We prove that every sufficiently large \(n\)-vertex \(k\)-uniform hypergraph
		with at least \(a n^h\) copies of \(H\) contains a blow-up of \(H\) with at
		least \(b(\log n)^{1/(k-1)}\) vertices in each part, where \(b=b(H,a)>0\).
		The exponent of $\log n$ is best possible. 
		This answers a question of R\"odl and Schacht and proves a conjecture of
		Bollob\'as, Erd\H{o}s and Simonovits from 1976.
		
		For fixed \(q\ge2\) and \(0<\alpha<1/2\), it follows that every
		\(q\)-coloring of the \(k\)-subsets of \([n]\) contains a set of order
		\((\log n)^{1/(k-1)}\) on which one color has density at least \(1-\alpha\).
For triples this recovers the sharp bound of Conlon, Fox and Sudakov. For \(k\ge4\), this improves their general exponent to the optimal value
\(1/(k-1)\) and answers a 1990 question of Erd\H{o}s for which he
offered a \$500 prize. 
	\end{abstract}
	
	\section{Introduction}
	
	Let \(H\) be a \(k\)-uniform hypergraph, or \(k\)-graph, on vertex set
	\([h]=\{1,\ldots,h\}\). The \emph{blow-up} \(H[t]\) is obtained by replacing
	each vertex \(i\) by a set \(U_i\) of \(t\) vertices and each edge
	\(\{i_1,\ldots,i_k\}\) by all \(t^k\) edges with one vertex in each of
	\(U_{i_1},\ldots,U_{i_k}\). The sets \(U_1,\ldots,U_h\) are disjoint.
	Write \(\cK\) for the complete \(k\)-graph on \(h\) vertices. We count
	copies as subhypergraphs, so that each copy of \(H\) gives
	\(|\aut(H)|\) injective edge-preserving maps. All logarithms are to base two.
	
	The Erd\H{o}s--Stone--Simonovits theorem, the fundamental theorem of
	extremal graph theory, can be stated in terms of blow-ups. Erd\H{o}s
	and Stone~\cite{ES46} proved that for fixed \(r\ge2\), \(t\ge1\) and
	\(0<\eps<1/(r-1)\), every sufficiently large graph of edge density at least
	\(1-1/(r-1)+\eps\) contains \(K_r[t]\). Every \(r\)-chromatic graph
	\(F\) is contained in some \(K_r[t]\), and the balanced complete
	\((r-1)\)-partite graph avoids \(F\). These observations give the Erd\H{o}s--Simonovits
	formula~\cite{ES66} for the Tur\'an density of \(F\), namely \(1-1/(r-1)\).
	 In 1973 Bollob\'as and Erd\H{o}s~\cite{BE}
	showed that for fixed \(r\) and \(\eps\) the correct order of magnitude of $t$ that we can guarantee is \(\log n\).
	
	This blow-up structure arises in other parts of mathematics as well, for example,  in additive combinatorics. For fixed \(k,t\)
	and \(\eta>0\), the multidimensional Szemer\'edi theorem of Furstenberg
	and Katznelson~\cite{FK} implies that every subset of \([N]^k\) of size
	at least \(\eta N^k\) contains a grid \(A_1\times\cdots\times A_k\)
	when \(N\) is sufficiently large. Here each \(A_i\) is a \(t\)-term
	arithmetic progression and all have the same positive common
	difference. To see why this is a blow-up, write \(S\) for the given subset of \([N]^k\).
	Form a \(k\)-partite \(k\)-graph with vertex classes
	\(V_i=\{i\}\times[N]\), in which each point \((x_1,\ldots,x_k)\in S\)
	gives the edge \(\{(1,x_1),\ldots,(k,x_k)\}\).
	Then \(A_1\times\cdots\times A_k\subseteq S\) precisely when every
	choice of one vertex from each of the sets \(\{i\}\times A_i\) is an
	edge. These sets are therefore the classes of a blow-up of a single
	\(k\)-edge.

	A major open question in extremal combinatorics is the following: If an \(n\)-vertex \(k\)-graph contains at least \(a n^h\) copies of \(H\), how large a blow-up of \(H\) must it contain? For graphs,
	Nikiforov~\cite{NikClique,NikGeneral} proved that all the parts can have
	order \(\log n\). Fox, Wigderson and Zhou~\cite{FWZ} recently improved the
	dependence on \(a\). For a single \(k\)-edge, the bound
	\(t=\Omega((\log n)^{1/(k-1)})\) is the dense case of the
	K\H{o}v\'ari--S\'os--Tur\'an theorem~\cite{KST} for \(k=2\) and a
	theorem of Erd\H{o}s from 1964~\cite{Erdos64} for \(k\ge3\).
	One can also apply the theorem of Erd\H{o}s to the \(h\)-graph whose edges are the
	vertex sets of copies of \(\cK\). This gives exponent \(1/(h-1)\).
	R\"odl and Schacht~\cite[Problem~2]{RS} asked whether the exponent
	\(1/(k-1)\) holds for every fixed \(H\). We answer this positively below.
	
	\begin{theorem}\label{thm:cliques}
		Fix integers \(h\ge k\ge2\) and a real \(a>0\). There are \(b>0\) and
		\(n_0\) such that every \(n\)-vertex \(k\)-graph with \(n\ge n_0\) and at
		least \(a n^h\) copies of \(\cK\) contains \(\cK[t]\) for some
		\(t\ge b(\log n)^{1/(k-1)}\).
	\end{theorem}
		As R\"odl and Schacht observed, the clique case implies the general case.
	We give the reduction in Section~\ref{sec:proof}.
	
	\begin{corollary}\label{cor:general}
		Let \(k\ge2\), let \(H\) be a fixed \(k\)-graph on \(h\ge1\) vertices and
		let \(a>0\). There are \(b=b(H,a)>0\) and \(n_0=n_0(H,a)\) such that every
		\(n\)-vertex \(k\)-graph with \(n\ge n_0\) and at least \(a n^h\) copies
		of \(H\) contains \(H[t]\) for some \(t\ge b(\log n)^{1/(k-1)}\).
	\end{corollary}
	
	In 1976 Bollob\'as, Erd\H{o}s and Simonovits~\cite{BES} conjectured that
	an edge density exceeding the Tur\'an density of \(\cK\) by a fixed
	positive amount forces a blow-up of \(\cK\) with parts of order
	\((\log n)^{1/(k-1)}\). By supersaturation, Theorem~\ref{thm:cliques}
	proves their conjecture. Section~\ref{sec:applications} gives the
	corresponding statement for every \(H\) with an edge.
	
	In 1990 Erd\H{o}s~\cite[pp.~21--22]{Erdos90} offered \$500 for
	determining how the order of the following two-color Ramsey function
	varies with \(\alpha\). For \(0\le\alpha<1/2\), he defined
	\(\F(n,\alpha)\) to be the least \(s\) for which there is a red-blue
	coloring of the \(k\)-subsets of \([n]\) such that, on every set \(S\) of
	at least \(s\) vertices, each color occurs on more than
	\(\alpha\binom{|S|}{k}\) of its \(k\)-subsets. At \(\alpha=0\) this is
	the inverse diagonal Ramsey function. For \(\alpha\) sufficiently close
	to \(1/2\), the order \((\log n)^{1/(k-1)}\) was already
	known~\cite{Erdos90}. Erd\H{o}s asked whether the change in the order of magnitude of
	\(\F(n,\alpha)\)
	as \(\alpha\) increases is continuous or occurs in jumps.\footnote{See also
		\href{https://www.erdosproblems.com/161}{Problem~161} in the Erd\H{o}s
		Problems database.}
	He suggested that the only jump might be at zero and singled out
	\(k>3\) as the range in which this might fail. Corollary~\ref{cor:jump} below shows that zero is indeed the only jump for every \(k\ge4\).
	
	In 2010 Conlon, Fox and Sudakov~\cite{CFS10} proved
	that for every \(k\ge2\) and every fixed \(0<\alpha<1/2\) there are
	\(c,\gamma>0\), depending on \(k\) and \(\alpha\), such that
	\(\F(n,\alpha)\ge c(\log n)^\gamma\). For \(k\ge4\), the
	Erd\H{o}s--Hajnal stepping-up bound gives \(\F(n,0)=O_k(\log\log n)\)
	(see~\cite[pp.~423--424]{CFS}). Together these bounds already show
	that the order of \(\F(n,\alpha)\) jumps at zero. In 2011 Conlon, Fox and
	Sudakov~\cite{CFS} proved that
	\(F_2^{(3)}(n,\alpha)=\Theta_\alpha(\sqrt{\log n})\) for every fixed
	\(0<\alpha<1/2\) and asked whether their result extends to \(k\ge4\). The proof of Corollary~\ref{cor:jump} gives this extension for every
	fixed number of colors. 
	\begin{corollary}\label{cor:jump}
		For each \(k\ge2\) and each fixed \(0<\alpha<1/2\),
		\[
		\F(n,\alpha)=\Theta_{k,\alpha}\bigl((\log n)^{1/(k-1)}\bigr).
		\]
		For \(k\ge4\), the value $\alpha=0$ is the only jump in the asymptotic order of
		\(\F(n,\alpha)\) on \(0\le\alpha<1/2\).
	\end{corollary}
	
	For graphs the order of magnitude of \(F_2^{(2)}(n,\alpha)\) is \(\log n\) even at \(\alpha=0\)~\cite{Erdos90}.
	For a \(k\)-graph \(H\), let \(r_k(H;q)\) be the least \(N\) such that
	every \(q\)-coloring of the \(k\)-subsets of \([N]\) contains a
	monochromatic copy of \(H\). Erd\H{o}s remarked that a bound
	\(F_2^{(3)}(n,\alpha)\ge c(\alpha)(\log n)^\gamma\) with \(\gamma>0\)
	independent of \(\alpha\) would make him doubt the conjectured
	double-exponential growth of \(r_3(K_s^{(3)};2)\)
	(see~\cite[Section~6.2]{CFS10}).
	Whether there is a jump at zero for triples remains open. Such a jump
	means that \(F_2^{(3)}(n,0)=o(\sqrt{\log n})\), or equivalently that
	\(\log r_3(K_s^{(3)};2)/s^2\to\infty\) as \(s\to\infty\).
	\medskip
	
	\noindent
	{\bf Notation and organization.}  
	Let \(V_1,\ldots,V_h\) be pairwise disjoint vertex classes.
	Following~\cite{ADLRY}, we call a product
	\(U_1\times\cdots\times U_h\), with \(U_i\subseteq V_i\), a
	\emph{cylinder}. If each \(U_i\) has size \(t\), we call it a
	\emph{\(t\)-cylinder}. Its \(t^h\) transversals are the sets obtained
	by choosing one vertex from each \(U_i\).
	
	In the proofs we count ordered lists of vertices, allowing repetitions
	within each list. For lists of length \(t\) there are \(t^h\) choices
	of positions, and several choices may give the same transversal.
	When the entries within each list are distinct, the lists specify a
	\(t\)-cylinder. We first count the choices of lists for which every
	transversal has the required property and then extract a cylinder
	with distinct vertices.
	
	Section~\ref{sec:warmups} gives the proof first for a single edge (to explain the notation used in the proof) and
	then for a triangle. We give the general proof in Section~\ref{sec:proof}.
	Section~\ref{sec:applications} gives the short proofs of the conjecture of
	Bollob\'as, Erd\H{o}s and Simonovits and of Corollary~\ref{cor:jump}.
	Appendix~\ref{app:tetrahedron} gives the calculation just in the case of the tetrahedron; this could be read immediately after Section~\ref{sec:warmups} to see the crucial inequalities involved, and how they combine in the proof.

	\section{A single edge and a triangle}\label{sec:warmups}
	
	We use the partition argument of Erd\H{o}s and Kleitman~\cite{EK}
	to pass to equal parts at the expense of only a constant multiplicative factor. If an \(n\)-vertex \(k\)-graph \(G\) contains
	\(A\) labeled embeddings of a \(k\)-graph \(H\) on \([h]\), where
	\(n\ge h\), then for \(N=\lfloor n/h\rfloor\) some disjoint \(N\)-sets
	\(V_1,\ldots,V_h\) satisfy
	\begin{equation}\label{eq:partition}
		\frac{\#\{\phi:H\hookrightarrow G \mid \phi(i)\in V_i
			\text{ for every }i\}}{N^h}
		\ge \frac{A}{h!{n \choose h}}.
	\end{equation}
	Indeed, over all ordered choices of the parts, every ordered \(h\)-tuple
	of distinct vertices occurs equally often as a transversal.
	
	\subsection{A single edge}\label{sec:edge}
	We give the proof for a single edge to introduce the notation used
	later. The common-neighborhood count is the one used in the
	K\H{o}v\'ari--S\'os--Tur\'an argument for graphs and in the proof of
	Erd\H{o}s for \(r>2\).
	
	Let \(G\) be an \(r\)-partite \(r\)-graph with parts \(V_1,\ldots,V_r\)
	of size \(N\) and \(pN^r\) edges, where \(r\ge2\) and \(p\ge p_0>0\).
	We seek \(t\) distinct vertices in each part such that all \(t^r\)
	transversals are edges. Fix \(2\le t\le N\). A choice \(X\) consists
	of an ordered list \((x_{i,1},\ldots,x_{i,t})\in V_i^t\) for each
	\(i\in[r]\). Its \emph{edge pattern} is
	\[
	J(X)=\{(a_1,\ldots,a_r)\in[t]^r:
	\{x_{1,a_1},\ldots,x_{r,a_r}\}\in E(G)\}.
	\]
	For a fixed set \(J\subseteq[t]^r\) of positions, define its density by
	\[
	\rho(J)=\frac{|\{X:J\subseteq J(X)\}|}{N^{rt}}.
	\]
	Thus \(J\) prescribes which transversals must be edges. The vertices
	occupying those positions vary with \(X\).
	Although \(\rho(J)\) is defined for every \(J\subseteq[t]^r\), for a
	single edge it suffices to compare the densities along the chain
	\[
	J_j=[t]^j\times\{1\}^{r-j},\qquad 0\le j\le r.
	\]
	At stage \(j\) we count every choice \(X\) for which all transversals
	in \(J_j\) are edges. 
	In \(J_j\) the first \(j\) coordinates range over all \(t\) positions,
	and each remaining coordinate is fixed at position \(1\). In particular,
	\(\rho(J_0)=p\), since every transversal occurs equally often in the
	first positions. We call \(X\) \emph{full} if \(J(X)=[t]^r\).
	
For \(r=2\), this chain consists of the three patterns	
	\[
	J_0=\{(1,1)\},\qquad J_1=[t]\times\{1\},\qquad J_2=[t]^2.
	\]
	The first requires the edge \(x_1y_1\), the second requires every
	\(x_u\) to be adjacent to \(y_1\), and the third requires every
	\(x_u\) to be adjacent to every \(y_v\). Write \(d(y)\) for the degree
	of \(y\in V_2\) and \(\Gamma(A)\subseteq V_2\) for the common
	neighborhood of \(A\subseteq V_1\). Then
	\begin{equation}\label{eq:stars}
		N^{t+1}\rho(J_1)
		=\sum_{y\in V_2}d(y)^t
		=\sum_{(x_1,\ldots,x_t)\in V_1^t}
		|\Gamma(\{x_1,\ldots,x_t\})|.
	\end{equation}
	Indeed, there are \(N^{t+1}\) choices of the first list and the center
	\(y_1\), and the other \(t-1\) entries of the second list are free.
	Counting stars with a set of \(t\) distinct leaves gives the usual
	K\H{o}v\'ari--S\'os--Tur\'an identity~\cite{KST},
	\(\sum_y\binom{d(y)}t=\sum_{A\in\binom{V_1}{t}}|\Gamma(A)|\).
	
	Regard \(J\subseteq[t]^r\) as an \(r\)-partite \(r\)-graph with a
	separate copy of \([t]\) in each class. The \emph{link} of position
	\(u\) in class \(i\) consists of the tuples obtained by deleting
	coordinate \(i\) from each tuple in \(J\) whose \(i\)th entry is \(u\).
	For instance, in the star pattern \(J_1\) for graphs, position \(1\)
	in the second class has link \([t]\), and every other position in
	that class has empty link. Giving each position the link \([t]\)
	yields \(J_2\). In general, for \(0\le j<r\), copying the link of
	position \(1\) in class \(j+1\) of \(J_j\) to every position in that
	class gives \(J_{j+1}\).
	
	Fix \(j\) in this range and let \(Y=(Y_i)_{i\in[r]\setminus\{j+1\}}\),
	where \(Y_i=(x_{i,1},\ldots,x_{i,t})\in V_i^t\). Thus \(Y\)
	specifies all lists except the one in \(V_{j+1}\). For
	\(\boldsymbol a=(a_1,\ldots,a_j)\in[t]^j\), put
	\[
	A_Y(\boldsymbol a)=\{x_{i,a_i}:1\le i\le j\}
	\cup\{x_{i,1}:j+2\le i\le r\}.
	\]
	This \((r-1)\)-set uses one chosen vertex from each of the first \(j\)
	lists and the first vertex of each of the last \(r-j-1\) lists. Let
	\[
	N(Y)=\{v\in V_{j+1}:A_Y(\boldsymbol a)\cup\{v\}\in E(G)
	\text{ for every }\boldsymbol a\in[t]^j\}
	\]
	and put \(d(Y)=|N(Y)|\). When \(j=0\), there is one empty tuple in
	\([t]^0\) and hence just one \((r-1)\)-set to complete.
	
	\medskip
	\noindent\textbf{Claim.}
	With both sums taken over the \(N^{(r-1)t}\) choices of \(Y\), we have
	\[
	\rho(J_j)=\frac{1}{N^{(r-1)t}}\sum_Y\frac{d(Y)}{N},
	\qquad
	\rho(J_{j+1})=\frac{1}{N^{(r-1)t}}\sum_Y
	\left(\frac{d(Y)}{N}\right)^t .
	\]
	\begin{proof}[Proof of the claim]
		Fix \(Y\). For \(J_j\), only the first entry in the missing list must
		belong to \(N(Y)\), giving \(d(Y)N^{t-1}\) choices. For \(J_{j+1}\)
		all \(t\) entries must belong to \(N(Y)\), giving \(d(Y)^t\) choices.
		Divide each count by the \(N^t\) completions of \(Y\) and average over
		\(Y\).
	\end{proof}
	By convexity of \(z^t\) on \([0,\infty)\) the claim gives
	\begin{equation}\label{eq:edge-copy}
		\rho(J_{j+1})\ge \rho(J_j)^t,\qquad
		\rho([t]^r)=\rho(J_r)\ge p^{\,t^r}.
	\end{equation}
	The second inequality follows by iterating the first from \(\rho(J_0)=p\).
	
	\medskip
	\noindent\textbf{Passing to distinct vertices.}
	The following count works for any prescribed family of transversals
	in \(V_1\times\cdots\times V_r\). Call a choice full when all its
	transversals belong to the family and write \(\rho([t]^r)\) for the
	fraction of full choices.
	Fix one class and fix the lists \(Y\) in all other classes. Let \(N(Y)\)
	be the set of vertices of the missing class that complete every
	transversal of \(Y\) to a member of the prescribed family, and put
	\(d(Y)=|N(Y)|\). There are \(d(Y)^t\) completions giving a full \(X\).
	If two specified positions in the missing list must also be equal,
	there are \(d(Y)^{t-1}\) such completions. The fraction of all choices
	\(X\) that are full and have this specified equality is consequently
	\[
	\frac{1}{N}\,\frac{1}{N^{(r-1)t}}\sum_Y
	\left(\frac{d(Y)}{N}\right)^{t-1}
	\le \frac{1}{N}
	\left(\frac{1}{N^{(r-1)t}}\sum_Y
	\left(\frac{d(Y)}{N}\right)^t\right)^{(t-1)/t}
	=\frac{\rho([t]^r)^{1-1/t}}{N},
	\]
	by concavity of \(z^{(t-1)/t}\). There are \(r\binom t2\) choices of a
	class and a pair of positions. If
	\begin{equation}\label{eq:distinct}
		N\,\rho([t]^r)^{1/t}>r\binom t2 ,
	\end{equation}
	the fraction of all \(X\) that are full and have a repeated entry
	is strictly less than \(\rho([t]^r)\). Hence some full choice has \(t\) distinct
	vertices in each class.
	
	For the single edge, \eqref{eq:edge-copy} gives
	\(N\rho([t]^r)^{1/t}\ge Np^{t^{r-1}}\).
	For sufficiently large \(N\), take
	\[
	t=\left\lfloor
	\left(\frac{\log N}{2\log(r/p_0)}\right)^{1/(r-1)}
	\right\rfloor .
	\]
	Then \(Np^{t^{r-1}}\ge\sqrt N>r\binom t2\), so \eqref{eq:distinct}
	gives the required \(t\)-cylinder. For \(r=2\) this gives
	\(K_{t,t}\) with \(t=\lfloor\log N/(2\log(2/p_0))\rfloor\).
	For a general \(n\)-vertex \(r\)-graph with at least \(a n^r\) edges,
	\eqref{eq:partition} gives \(N=\lfloor n/r\rfloor\) and \(p_0=r!a\).
	Thus \(t=\Omega_{r,a}((\log n)^{1/(r-1)})\), recovering the dense
	case of Erd\H{o}s's theorem~\cite{Erdos64} for \(r\ge3\).
	
	If we regard the triangles of a tripartite graph as edges of a
	\(3\)-graph, this argument gives parts of order \(\sqrt{\log N}\).
	To obtain parts of order \(\log N\), the next subsection uses the fact
	that the three bipartite graphs on the lists determine the triangle pattern.
	
	\subsection{Triangles}\label{sec:triangles}
	We now give the proof for \(H=K_3\), a result already known from the work of
	Nikiforov. Let \(G\) be a tripartite graph with parts \(V_1,V_2,V_3\)
	of size \(N\), and suppose it has \(pN^3\) triangles, where \(p\ge p_0>0\).
	We prove that, for \(N\) sufficiently large in terms of \(p_0\), \(G\)
	contains a \(K_{t,t,t}\), with one class in each \(V_i\), for
	\begin{equation}\label{eq:tri-size}
		t=\left\lfloor\frac{p_0}{3}(\log N-2\log\log N)\right\rfloor .
	\end{equation}
	For \(a n^3\) triangles, \eqref{eq:partition} gives \(p_0=6a\) and
	\(N=\lfloor n/3\rfloor\), so \(t=(2a-o(1))\log n\). This gives the
	logarithmic order in Nikiforov's theorem~\cite{NikClique}.
	
	We use the following consequence of the arithmetic--geometric mean
	inequality (AM--GM). For nonnegative functions \(g_1,\ldots,g_t\)
	on a nonempty finite set \(\mathcal Y\),
	\begin{equation}\label{eq:amgm-product}
		\frac{1}{|\mathcal Y|}\sum_{Y\in\mathcal Y}\prod_{u=1}^t g_u(Y)
		\le \prod_{u=1}^t
		\left(\frac{1}{|\mathcal Y|}\sum_{Y\in\mathcal Y} g_u(Y)^t\right)^{1/t}.
	\end{equation}
	If a factor on the right is zero, then the corresponding function
	vanishes identically. Otherwise normalize each \(g_u\) so that its
	\(t\)th power has average one. For each \(Y\), AM--GM gives
	\(\prod_u g_u(Y)\le t^{-1}\sum_u g_u(Y)^t\). Averaging over \(Y\)
	gives the normalized inequality, and rescaling gives~\eqref{eq:amgm-product}.
	
	\begin{proof}[Proof for triangles]
		Fix \(2\le t\le N\). Choose three ordered lists
		\((x_1,\ldots,x_t)\in V_1^t\), \((y_1,\ldots,y_t)\in V_2^t\) and
		\((z_1,\ldots,z_t)\in V_3^t\), allowing repetitions. Write \(\cX\)
		for the family of the \(N^{3t}\) choices \(X\).
		A triple of positions \((u,v,w)\in[t]^3\) selects the vertices
		\(x_u,y_v,z_w\). Define the \emph{triangle pattern} of \(X\) by
		\[
		J(X)=\{(u,v,w)\in[t]^3:x_u y_v z_w\text{ is a triangle of }G\}.
		\]
		The edges between the three pairs of lists determine \(J(X)\), since
		\((u,v,w)\in J(X)\) precisely when \(x_u y_v,x_u z_w,y_v z_w\)
		are all edges. Thus \(3t^2\) binary entries determine the pattern,
		and the number \(Q\) of patterns that occur satisfies \(Q\le 2^{3t^2}\).
		
		For a fixed \(J\subseteq[t]^3\), its density is
		\[
		\rho(J)=\frac{|\{X\in\cX:J\subseteq J(X)\}|}{N^{3t}}.
		\]
		Only the vertices occupying the positions vary with \(X\). Triples of
		positions outside \(J\) are unrestricted. The fraction of full choices
		is \(\rho([t]^3)\). In a full choice, the sets of vertices appearing
		in the three lists span a complete tripartite graph. If every list has
		\(t\) distinct entries, this graph is a \(K_{t,t,t}\).
		For example, if \(J=\{(1,1,1),(1,2,1)\}\), the two required triangles
		are \(x_1y_1z_1\) and \(x_1y_2z_1\). Both use the edge \(x_1z_1\),
		as in Figure~\ref{fig:triangles}.
		
		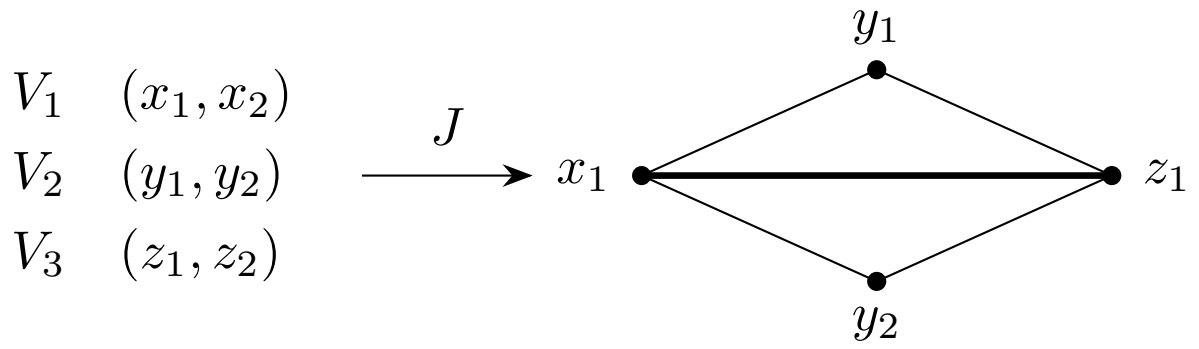
\begin{figure}[htbp]
			\centering
			\begin{tikzpicture}[x=1cm,y=1cm,font=\small]
				\node at (-2.52,0) {$\begin{array}{c@{\quad}l}
						V_1 & (x_1,x_2)\\[3pt]
						V_2 & (y_1,y_2)\\[3pt]
						V_3 & (z_1,z_2)
					\end{array}$};
				\draw[-{Stealth[length=2mm]}] (-1.13,0) -- (0.00,0)
				node[midway,above=2pt] {\(J\)};
				\coordinate (x) at (.72,0);
				\coordinate (z) at (3.83,0);
				\coordinate (p) at (2.275,.70);
				\coordinate (q) at (2.275,-.70);
				\draw (x)--(p)--(z)--(q)--cycle;
				\draw[line width=1.25pt] (x)--(z);
				\foreach \point in {x,z,p,q}
				\fill (\point) circle[radius=1.8pt];
				\node[left=2pt] at (x) {\(x_1\)};
				\node[right=2pt] at (z) {\(z_1\)};
				\node[above=1pt] at (p) {\(y_1\)};
				\node[below=1pt] at (q) {\(y_2\)};
			\end{tikzpicture}
			\caption{For \(J=\{(1,1,1),(1,2,1)\}\), the requirements give
				the five edges shown when \(y_1\ne y_2\). The two triangles
				share the edge \(x_1z_1\).}
			\label{fig:triangles}
		\end{figure}
		
		For \(u\in[t]\), the \emph{link} of position \(u\) in the first list is
		\[
		J_u=\{(v,w)\in[t]^2:(u,v,w)\in J\}.
		\]
		Thus \(J_u\) specifies the pairs \(y_vz_w\) that \(x_u\) must complete
		to a triangle. Imposing the requirements of position \(u\) at every
		position of the first list gives the pattern
		\[
		S_{1,u}J=[t]\times J_u=\{(a,v,w)\in[t]^3:(v,w)\in J_u\}.
		\]
		For example, if \(J_u=\{(1,2)\}\), then \(S_{1,u}J\) requires
		\(x_a y_1z_2\) to be a triangle for every \(a\in[t]\).
		The operations in the other two lists are
		\[
		S_{2,v}J=\{(a,b,c)\in[t]^3:(a,v,c)\in J\},\qquad
		S_{3,w}J=\{(a,b,c)\in[t]^3:(a,b,w)\in J\}.
		\]
		
		Fix the second and third lists \(Y=((y_1,\ldots,y_t),(z_1,\ldots,z_t))\).
		For each \(u\in[t]\), let
		\[
		N_u(Y)=\bigl\{x\in V_1:
		\{xy_v,xz_w,y_vz_w\}\subseteq E(G)
		\text{ for every }(v,w)\in J_u\bigr\},
		\qquad d_u(Y)=|N_u(Y)|.
		\]
		Thus a single vertex \(x\in N_u(Y)\) completes every pair prescribed
		at position \(u\). If a prescribed pair \(y_vz_w\) is missing,
		\(N_u(Y)=\varnothing\). If \(J_u=\varnothing\), then \(N_u(Y)=V_1\).
		
		For \(J\) we choose \(x_1\) from \(N_1(Y)\), \(x_2\) from \(N_2(Y)\),
		and so on. The possible first lists form the product
		\(N_1(Y)\times\cdots\times N_t(Y)\). For \(S_{1,u}J\), fix \(u\)
		and choose all \(t\) entries from the one set \(N_u(Y)\). The possible
		first lists now form \(N_u(Y)^t\), one such product for each choice of
		\(u\). Both counts are exact even when the sets overlap, since the
		entries may repeat. Summing over the \(N^{2t}\) choices of \(Y\) gives
		\begin{equation}\label{eq:triangle-counts}
			\begin{aligned}
				\rho(J)&=\frac{1}{N^{2t}}\sum_Y\prod_{u=1}^t\frac{d_u(Y)}{N},\\
				\rho(S_{1,u}J)&=\frac{1}{N^{2t}}\sum_Y
				\left(\frac{d_u(Y)}{N}\right)^t .
			\end{aligned}
		\end{equation}
		
		The second identity also occurs in dependent random choice~\cite{FS}.
		Join \(x\in V_1\) to
		\((y,z)\in V_2\times V_3\) when \(xyz\) is a triangle. Then \(N_u(Y)\)
		is the common neighborhood of the pairs \((y_v,z_w)\) prescribed by
		\(J_u\), so exactly \(d_u(Y)^t\) ordered lists have all their entries
		in this common neighborhood.
		
		Apply~\eqref{eq:amgm-product} to the first line
		of~\eqref{eq:triangle-counts}, with \(g_u(Y)=d_u(Y)/N\).
		The second line identifies the average of \(g_u(Y)^t\) as
		\(\rho(S_{1,u}J)\). The same count works with any one of the three
		lists left to choose. Thus, for \(i=1,2,3\),
		\begin{equation}\label{eq:copy}
			\rho(J)\le\prod_{u=1}^t \rho(S_{i,u}J)^{1/t}.
		\end{equation}
		For \(u,v\in[t]\), put \(J_{u,v}=\{w\in[t]:(u,v,w)\in J\}\).
		After copying at \(u\) in the first list and at \(v\) in the second,
		we have \(S_{2,v}(S_{1,u}J)=[t]^2\times J_{u,v}\). Indeed,
		\((a,b,w)\) belongs to this pattern precisely when \((u,v,w)\in J\).
		A further copy at \(w\) in the third list gives \([t]^3\) if
		\(w\in J_{u,v}\), and the empty pattern otherwise. Applying
		\eqref{eq:copy} in the three lists in turn therefore gives
		\begin{align*}
			\rho(J)
			&\le \prod_{u\in[t]}\rho([t]\times J_u)^{1/t}\\
			&\le \prod_{u,v\in[t]}\rho([t]^2\times J_{u,v})^{1/t^2}\\
			&\le \prod_{(u,v,w)\in J}\rho([t]^3)^{1/t^3}
			\prod_{(u,v,w)\in[t]^3\setminus J}\rho(\varnothing)^{1/t^3}.
		\end{align*}
		Since \(\rho(\varnothing)=1\), we obtain
		\begin{equation}\label{eq:tri-box}
			\rho(J)\le\rho([t]^3)^{|J|/t^3}.
		\end{equation}
		For \(J=\{(1,1,1)\}\) we have \(\rho(J)=p\), so~\eqref{eq:tri-box}
		already gives \(\rho([t]^3)\ge p^{t^3}\). We now use the fact that
		at most \(2^{3t^2}\) triangle patterns occur to obtain the bound
		\(\rho([t]^3)\ge 2^{-3t^2/p}\). For each \(J\subseteq[t]^3\), let
		\(n_J\) be the number of choices \(X\in\cX\) for which \(J(X)=J\).
	Every such \(X\) contributes to the numerator defining \(\rho(J)\), so \(n_J\le|\cX|\rho(J)\). Grouping by the exact
		pattern, we obtain
		\begin{equation}\label{eq:tri-exact}
			\frac{1}{|\cX|}\sum_{X\in\cX}\frac{1}{\rho(J(X))}
			=\sum_{J:\,n_J>0}\frac{n_J}{|\cX|\rho(J)}
			\le Q .
		\end{equation}
			Since $\rho(J(X))>0$ for all $X$, applying AM--GM  gives
		\[
		\left(\prod_{X\in\cX}\frac{1}{\rho(J(X))}\right)^{1/|\cX|}
		\le \frac{1}{|\cX|}\sum_{X\in\cX}\frac{1}{\rho(J(X))}
		\le Q.
		\]
		Taking reciprocals in the preceding display and then applying
		\eqref{eq:tri-box} to each factor \(\rho(J(X))\) gives
			\[
		\frac{1}{Q}
		\le \left(\prod_{X\in\cX}\rho(J(X))\right)^{1/|\cX|}
		\le \left(\prod_{X\in\cX}
		\rho([t]^3)^{|J(X)|/t^3}\right)^{1/|\cX|}
		= \rho([t]^3)^{\frac{1}{t^3|\cX|}\sum_{X\in\cX}|J(X)|}.
		\]
		To evaluate the exponent, fix a triple of positions \((u,v,w)\in[t]^3\).
		There are \(pN^3\) ways to place a triangle at these positions and
		\(N^{3(t-1)}\) ways to fill the remaining positions. Thus
		\((u,v,w)\in J(X)\) for exactly \(pN^{3t}=p|\cX|\) choices \(X\).
		Summing over the \(t^3\) triples of positions gives
		\(\sum_{X\in\cX}|J(X)|=p\,t^3|\cX|\). The exponent is therefore \(p\),
		so \(1/Q\le \rho([t]^3)^p\). Since \(Q\le 2^{3t^2}\), we obtain
		\[
		\rho([t]^3)\ge Q^{-1/p}\ge 2^{-3t^2/p}.
		\]
		
		Finally, take \(t\) as in~\eqref{eq:tri-size}. Since \(p\ge p_0\),
		we have \(3t/p\le \log N-2\log\log N\). Also,
		\(t\le (p_0/3)\log N\le (\log N)/3\). Hence
		\[
		N\rho([t]^3)^{1/t}
		\ge N2^{-3t/p}
		\ge (\log N)^2
		> 3\binom{t}{2}.
		\]
		By~\eqref{eq:distinct}, some full choice has \(t\) distinct vertices
		in each list. Every transversal of these three sets is a triangle,
		so the sets span a \(K_{t,t,t}\).
		\end{proof}

	\section{Proof of the blow-up theorem}\label{sec:proof}
In this section, we provide the formal proof of Theorem~\ref{thm:cliques}.  If the reader has  absorbed the proofs for an edge and triangle, then there are essentially no new mathematical ideas except to extend the notation from the previous section to the general setting. For a more streamlined presentation without the general definitions, the reader is advised to see the Appendix, where a proof is presented for tetrahedra.

	By~\eqref{eq:partition} we may work with \(h\) vertex classes
	\(V_1,\ldots,V_h\) of size \(N\). For a \(k\)-graph \(G\) on these
	classes, let \(T\subseteq V_1\times\cdots\times V_h\) consist of the
	tuples \((v_1,\ldots,v_h)\) that span a copy of \(\cK\). The parts
	fix the order of the vertices, so \(|T|\) counts the transversal
	\(h\)-cliques of \(G\). We seek \(t\)-sets \(U_i\subseteq V_i\)
	with \(U_1\times\cdots\times U_h\subseteq T\). This means that all
	\(t^h\) transversals are cliques, giving the desired blow-up \(\cK[t]\).
	
	We first extend the count from Section~\ref{sec:triangles} to an
	arbitrary family \(T\subseteq V_1\times\cdots\times V_h\) of density
	\(p=|T|/N^h>0\). Fix \(t\ge2\). Let
	\(\cX=V_1^t\times\cdots\times V_h^t\), so there are \(N^{ht}\)
	choices \(X\) of the \(h\) ordered lists. Write the list in
	\(V_i\) as \((x_{i,1},\ldots,x_{i,t})\). The pattern of \(X\) and the
	density of a fixed \(J\subseteq[t]^h\) are given by
	\[
	J(X)=\{\boldsymbol a\in[t]^h:(x_{1,a_1},\ldots,x_{h,a_h})\in T\},
	\qquad
	\rho(J)=\frac{|\{X\in\cX:J\subseteq J(X)\}|}{N^{ht}}.
	\]
We call \(X\) \emph{full} if \(J(X)=[t]^h\), that is, if every
choice of one entry from each of the \(h\) lists gives a member of \(T\).	Let \(Q\) be the number of distinct patterns \(J(X)\). Repeating
	the entries of any one member of \(T\) in their respective lists
	gives a full choice, so \(\rho([t]^h)>0\).
	
	\begin{lemma}\label{lem:cylinders}
		The density of the full pattern satisfies
		\begin{equation}\label{eq:criterion}
			\rho([t]^h)\ge Q^{-1/p}.
		\end{equation}
	\end{lemma}
	
	\begin{proof}
		Fix \(J\subseteq[t]^h\) and \(i\in[h]\). Copying the link of position
		\(u\) in class \(i\) to every position in that class gives the pattern
		\[
		S_{i,u}J=\{\boldsymbol a\in[t]^h:
		(a_1,\ldots,a_{i-1},u,a_{i+1},\ldots,a_h)\in J\}.
		\]
		Fix the \(h-1\) lists \(Y\) outside \(V_i\). Let \(N_u(Y)\subseteq V_i\)
		be the set of vertices \(v\) with the following property. For each
		\(\boldsymbol a\in J\) with \(a_i=u\), the tuple having \(v\) in
		coordinate \(i\) and \(x_{\ell,a_\ell}\) in each coordinate
		\(\ell\ne i\) belongs to \(T\). Put \(d_u(Y)=|N_u(Y)|\).
		
		As in~\eqref{eq:triangle-counts}, the completions for \(J\) form
		\(N_1(Y)\times\cdots\times N_t(Y)\), and the completions for
		\(S_{i,u}J\) form \(N_u(Y)^t\). Averaging the two counts over the
		\(N^{(h-1)t}\) choices of \(Y\) and using~\eqref{eq:amgm-product} gives
		\[
		\begin{aligned}
			\rho(J)
			&=\frac{1}{N^{(h-1)t}}\sum_Y\prod_{u=1}^t\frac{d_u(Y)}{N}\\
			&\le\prod_{u=1}^t
			\left(\frac{1}{N^{(h-1)t}}\sum_Y
			\left(\frac{d_u(Y)}{N}\right)^t\right)^{1/t}
			=\prod_{u=1}^t\rho(S_{i,u}J)^{1/t}.
		\end{aligned}
		\]
		Copy successively in classes \(1,\ldots,h\), applying this inequality
		to each factor at the next step. After \(\ell\) steps this gives
		\[
		\rho(J)\le\prod_{j_1,\ldots,j_\ell\in[t]}
		\rho\bigl(S_{\ell,j_\ell}\cdots S_{1,j_1}J\bigr)^{1/t^\ell}.
		\]
		A tuple \(\boldsymbol a\) belongs to the pattern in this product
		precisely when \((j_1,\ldots,j_\ell,a_{\ell+1},\ldots,a_h)\in J\).
		For \(\ell=h\) the pattern is full for the \(|J|\) choices of
		\((j_1,\ldots,j_h)\) in \(J\) and empty for the others. As
		\(\rho(\varnothing)=1\), this gives the analogue of~\eqref{eq:tri-box},
		\begin{equation}\label{eq:box}
			\rho(J)\le\rho([t]^h)^{|J|/t^h}.
		\end{equation}
			To apply~\eqref{eq:box} to \(J(X)\), we need two equations. First,
		write \(n_J=|\{X:J(X)=J\}|\). Since \(n_J\le|\cX|\rho(J)\),
		grouping by exact pattern as in~\eqref{eq:tri-exact} gives
		\[
		\frac{1}{|\cX|}\sum_{X\in\cX}\frac{1}{\rho(J(X))}
		=\sum_{J:\,n_J>0}\frac{n_J}{|\cX|\rho(J)}
		\le Q.
		\]
		Second, at any fixed \(\boldsymbol a\in[t]^h\), each member of \(T\)
		can occupy those positions in exactly \(N^{h(t-1)}\) choices of \(X\).
		Thus a fraction \(p\) of all choices have \(\boldsymbol a\in J(X)\).
		Summing over \(\boldsymbol a\) gives
		\[
		\frac{1}{|\cX|}\sum_{X\in\cX}\frac{|J(X)|}{t^h}=p .
		\]
		The first equation and AM--GM bound the geometric mean of the
		\(\rho(J(X))\) from below. The second identity gives the exponent when
		we take the geometric mean of~\eqref{eq:box}. Together they give
		\[
		Q^{-1}\le\left(\prod_{X\in\cX}\rho(J(X))\right)^{1/|\cX|}
		\le \rho([t]^h)^p .
		\]
		Taking the \(p\)th root proves~\eqref{eq:criterion}.
	\end{proof}
	
	\begin{proof}[Proof of Theorem~\ref{thm:cliques}]
		Return to the transversal \(h\)-cliques of \(G\), and suppose their
		density is \(p\ge p_0>0\). Put \(m=\binom hk\). There are \(m\) choices
		of \(k\) parts and \(t^k\) choices of one position in each. Record for
		each of these \(mt^k\) choices whether its vertices form an edge.
		A transversal is a clique precisely when all its \(m\) edges occur.
		Thus the entries determine \(J(X)\), giving
		\[
		Q\le 2^{mt^k},\qquad \rho([t]^h)\ge 2^{-mt^k/p}
		\]
		by Lemma~\ref{lem:cylinders}. A full choice with at least \(t\)
		distinct entries in every list gives a copy of \(\cK[t]\). Indeed,
		select \(t\) such vertices in each class. Any \(k\) selected vertices
		from distinct classes extend to a transversal clique and hence form
		an edge.
		
	Fix the part size \(t\) that we want in the theorem. The preceding estimate holds for every
	list length \(s\). Applying it with \(s=2t\) gives
	\[
	\rho([2t]^h)\ge 2^{-m(2t)^k/p}.
	\]
	It remains to show that among the choices counted by \(\rho([2t]^h)\),
	some choice has at least \(t\) distinct vertices in each list. If there is
		no such blow-up, then every full choice  has a list with at most \(t-1\) distinct entries. All entries of that list therefore lie in some \(t\)-subset of its class. Counting the possible classes and \(t\)-sets gives
		\[
		2^{-m(2t)^k/p}\le\rho([2t]^h)
		\le h\binom Nt(t/N)^{2t}\le h(et/N)^t.
		\]
		Taking \(t\)th roots gives \(N2^{-2^kmt^{k-1}/p}\le h^{1/t}et\).
		For \(t=\lfloor c(\log N)^{1/(k-1)}\rfloor\), with \(c>0\) chosen
		so that \(2^kmc^{k-1}\le p_0/2\), the left side is at least
		\(\sqrt N\), a contradiction for large \(N\).
		Thus we can choose a list of \(2t\) vertices from each \(V_i\) so that
		each list contains at least \(t\) distinct vertices and every selection
		of one entry from each of the \(h\) lists spans a copy of \(\cK\).
		For each \(i\), select \(t\) distinct vertices from the \(i\)th list
		and call the resulting set \(U_i\). Every transversal of \(U_1,\ldots,U_h\)
		then spans a clique, so these sets give a copy of \(\cK[t]\).

		Now suppose that the original \(n\)-vertex \(k\)-graph contains at
		least \(a n^h\) copies of \(\cK\). Each copy gives \(h!\) labeled
		embeddings, so by~\eqref{eq:partition} there are disjoint sets
		\(V_1,\ldots,V_h\), each of size \(N=\lfloor n/h\rfloor\), whose
		transversal clique density satisfies
		\[
		p\ge \frac{h!\,a n^h}{h!{n \choose h}}\ge h!\,a.
		\]
		Applying the partite bound with \(p_0=h!\,a\) gives a copy of
		\(\cK[t]\) with \(t=\lfloor c(\log N)^{1/(k-1)}\rfloor\), where
		\(c=c(h,k,a)>0\). Since \(\log N=\log n-O_h(1)\), we have
		\(t\ge b(\log n)^{1/(k-1)}\) for some \(b=b(h,k,a)>0\) and all
		sufficiently large \(n\). This proves the theorem.
		\end{proof}

	\begin{proof}[Proof of Corollary~\ref{cor:general}]
		Let \(G\) be \(h\)-partite with parts \(V_i\) of size \(N\) and
		transversal \(H\)-density at least \(p_0>0\). If \(H\) has no edges,
		take \(t=N\). Otherwise, for each \(I\in\binom{[h]}k\setminus E(H)\),
		add all crossing \(k\)-sets on the parts indexed by \(I\), and call
		the resulting \(k\)-graph \(G'\). A transversal is a copy of \(H\)
		in the prescribed order in \(G\) exactly when it is an \(h\)-clique
		in \(G'\). The partite form of Theorem~\ref{thm:cliques} therefore
		gives a \(\cK[t]\) in \(G'\) with its \(i\)-th class in \(V_i\).
		All edges of this blow-up prescribed by \(H\) already belong to
		\(G\), giving \(H[t]\). Now apply~\eqref{eq:partition} with
		\(p_0=|\aut(H)|a\) to obtain the corollary.
	\end{proof}
	
	\section{Extremal and Ramsey consequences}\label{sec:applications}
	
	\paragraph{The Tur\'an density.}
	Fix a \(k\)-graph \(H\) with \(h\) vertices, at least one edge and
	Tur\'an density \(\pi(H)\). For each \(\eps>0\), the supersaturation theorem
	of Erd\H{o}s and Simonovits~\cite[Theorem~\(1^*\)]{ES83} gives at least \(a(H,\eps)n^h\) copies of
	\(H\) in every sufficiently large \(n\)-vertex \(k\)-graph of density
	at least \(\pi(H)+\eps\), so Corollary~\ref{cor:general} gives \(H[t]\)
	with \(t\ge b(H,\eps)(\log n)^{1/(k-1)}\).

	\paragraph{Sharpness.}
	The exponent in Corollary~\ref{cor:general} is sharp for every
	\(h\)-vertex \(k\)-graph \(H\) with an edge, even for \(h\)-partite
	hosts. Take \(h\) parts of size \(N\) and include each crossing
	\(k\)-set independently with probability \(1/2\). The expected number
	of copies of \(H\) with vertex \(i\) in part \(i\) is
	\(2^{-|E(H)|}N^h\). Since there are at most \(N^h\) such copies,
	with probability bounded away from zero there are at least
	\(2^{-|E(H)|-1}N^h\) of them.
	Now let \(t=\lceil C_k(\log N)^{1/(k-1)}\rceil\). As in the
	probabilistic construction of Erd\H{o}s~\cite[Theorem~3]{Erdos64},
	the expected number of blow-ups of a single edge with parts of size
	\(t\) is at most \((hN)^{kt}2^{-t^k}=o(1)\) for \(C_k\) sufficiently
	large. Thus some host has \(\Omega_H(N^h)\) transversal copies of
	\(H\) and contains no such blow-up. Since \(H\) has an edge, this
	host contains no \(H[t]\).

	\paragraph{Monochromatic blow-ups.}
	For fixed \(q\ge2\) and \(h\ge k\ge2\), we have
	\begin{equation}\label{eq:ramsey}
		r_k(\cK[t];q)=2^{\Theta_{h,k,q}(t^{k-1})}.
	\end{equation}
	Conlon, Fox and Sudakov~\cite{CFS10} used Ramsey bounds for
	these blow-ups to find almost monochromatic sets. They proved the sharp
	exponent in~\eqref{eq:ramsey} for triples in~\cite{CFS}.
	For the upper bound, take \(R=r_k(\cK;q)\) disjoint parts of size
	\(N=\lfloor n/R\rfloor\) in a \(q\)-coloring of the \(k\)-sets of \([n]\).
	Every transversal of the \(R\) parts contains a monochromatic \(h\)-clique.
	There are \(q\binom Rh\) choices of its color and its \(h\) parts.
	Some choice occurs in at least \(N^R/(q\binom Rh)\) transversals.
	Each clique of that color on those parts extends to \(N^{R-h}\)
	transversals, so its color has transversal clique density at least
	\(1/(q\binom Rh)\) there. The partite form of Theorem~\ref{thm:cliques}
	gives a monochromatic \(\cK[t]\) with
	\(t\ge b(h,k,q)(\log n)^{1/(k-1)}\).
	The lower bound follows from the two-color construction of
	Erd\H{o}s~\cite{Erdos64}, since every \(\cK[t]\)
	contains a blow-up of a single \(k\)-edge with parts of size \(t\).
	
	\begin{proof}[Proof of Corollary~\ref{cor:jump}]
		More generally, fix \(q\ge2\) and color the \(k\)-sets of \([n]\)
		with \(q\) colors. Choose \(h\ge k\) with \(\binom k2/h\le\alpha\).
		By~\eqref{eq:ramsey}, there is a monochromatic \(\cK[t]\) with
		\(t\ge b(k,q,\alpha)(\log n)^{1/(k-1)}\). In the union of its \(h\)
		parts, the proportion of \(k\)-sets meeting \(k\) different parts is
		\[
		\frac{\binom hk t^k}{\binom{ht}k}
		\ \ge\ \prod_{j=0}^{k-1}\left(1-\frac{j}{h}\right)
		\ \ge\ 1-\frac{\binom k2}{h}
		\ \ge\ 1-\alpha .
		\]
		All these \(k\)-sets have the same color. In particular, when \(q=2\)
		the other color has density at most \(\alpha\) on this set. This gives
		the lower bound on \(\F(n,\alpha)\).
		
		For the upper bound, color each \(k\)-subset of \([n]\) independently
		red or blue with probability \(1/2\) each. Let
		\(s=\lceil C(\log n)^{1/(k-1)}\rceil\), put \(u=\binom sk\), and let
		\(Z\) be the number of \(s\)-sets on which some color has density at most
		\(\alpha\). For a fixed \(s\)-set \(S\), the number \(R_S\) of red edges
		has distribution \(\operatorname{Bin}(u,1/2)\). We use the lower-tail
		Chernoff bound
		\[
		\Pr\bigl(W\le(1-\xi)\mu\bigr)\le\exp(-\xi^2\mu/2)
		\qquad(0\le\xi\le1)
		\]
		for a binomial random variable \(W\) of mean \(\mu\). Taking
		\(\mu=u/2\) and \(\xi=1-2\alpha\) gives
		\[
		2^{-u}\sum_{j=0}^{\lfloor\alpha u\rfloor}\binom uj
		=\Pr(R_S\le\alpha u)
		\le \exp\!\left(-\frac{(1-2\alpha)^2u}{4}\right)
		=\exp\!\bigl(-(\tfrac12-\alpha)^2u\bigr).
		\]
		The two colors have the same distribution, and since \(\alpha<1/2\)
		their lower-tail events on \(S\) are disjoint. Linearity of
		expectation gives
		\[
		\mathbb E Z=2\binom ns\,\Pr(R_S\le\alpha u)
		\le 2\binom ns\exp\!\bigl(-(\tfrac12-\alpha)^2u\bigr).
		\]
		Since \(u/(s\log n)\to C^{k-1}/k!\), taking \(C=C(k,\alpha)\) large
		enough makes this less than \(1\) for all sufficiently large \(n\).
		Some coloring therefore has \(Z=0\). Both colors then have density
		greater than \(\alpha\) on every larger set \(W\) as well, since each
		\(k\)-subset of \(W\) lies in exactly \(\binom{|W|-k}{s-k}\) of its
		\(s\)-subsets.
		
		Finally, let \(k\ge4\). The classical stepping-up bound gives
		\(r_k(K_s^{(k)};2)\ge 2^{2^{c_k s}}\) for large \(s\)
		(see~\cite[pp.~423--424]{CFS}). Taking the inverse at \(\alpha=0\) gives
		\[
		\F(n,0)=O_k(\log\log n)=o\bigl((\log n)^{1/(k-1)}\bigr).
		\qedhere
		\]
	\end{proof}
	
		\section*{AI Disclosure}
	
	AI tools (GPT 6.2 Astra) generated the main proof ideas in this project.
	Through repeated discussions with the AI, I obtained simplifications of
	the proof. I then revised the proofs based on my understanding and
	simplified and changed them further. 
	AI also helped write and revise the manuscript and produced Figure~\ref{fig:triangles}. I have independently checked all proofs
	and take full responsibility for the content of this paper.

	\appendix
	\section{The tetrahedron}\label{app:tetrahedron}
	
Here we give a proof just for tetrahedra (without re-stating definitions), to explain to the interested reader what the crucial inequalities in the general proof are. 	Suppose the four parts have size \(N\) and at least \(pN^4\) transversals
are tetrahedra, where \(p>0\) is fixed. We show that for large \(N\)
there is a \(K_4^{(3)}[t]\) with \(t\) of order \(\sqrt{p\log N}\).	 Use the notation of Section~\ref{sec:proof} with
	\(s\) entries per list. 
	\begin{proof}
		Fix \(i\in[4]\). Counting the possible entries in list \(i\) with
		the other three lists fixed, and applying~\eqref{eq:amgm-product}, gives
		\[
		\rho(J)=\frac{1}{N^{3s}}\sum_Y\prod_{u=1}^s\frac{d_u(Y)}{N}
		\le\prod_{u=1}^s
		\left(\frac{1}{N^{3s}}\sum_Y
		\left(\frac{d_u(Y)}{N}\right)^s\right)^{1/s}
		=\prod_{u=1}^s\rho(S_{i,u}J)^{1/s}.
		\]
		Copying in classes \(1,2,3,4\) in turn gives
		\begin{align*}
			\rho(J)
			&\le\prod_{a\in[s]}\rho(S_{1,a}J)^{1/s}\\
			&\le\prod_{a,b\in[s]}\rho(S_{2,b}S_{1,a}J)^{1/s^2}\\
			&\le\prod_{a,b,c\in[s]}\rho(S_{3,c}S_{2,b}S_{1,a}J)^{1/s^3}\\
			&\le\prod_{a,b,c,d\in[s]}\rho(S_{4,d}S_{3,c}S_{2,b}S_{1,a}J)^{1/s^4}
			=\rho([s]^4)^{|J|/s^4}.
		\end{align*}
			Fix \((a,b,c,d)\in[s]^4\). After the first operation, a quadruple
		\((i,j,k,\ell)\) belongs to \(S_{1,a}J\) if and only if
		\((a,j,k,\ell)\in J\). Applying the same observation in the other
		three coordinates gives
		\[
		S_{4,d}S_{3,c}S_{2,b}S_{1,a}J
		=
		\begin{cases}
			[s]^4, & \text{if } (a,b,c,d)\in J,\\
			\varnothing, & \text{if } (a,b,c,d)\notin J.
		\end{cases}
		\]
		Thus each of the \(|J|\) quadruples in \(J\) contributes a factor
		\(\rho([s]^4)^{1/s^4}\) to the last product, and each of the remaining
		quadruples contributes \(\rho(\varnothing)^{1/s^4}=1\).
			There are \(4s^3\) ways to choose three of the four lists and one
		position in each. These edge indicators determine \(J(X)\), so the number \(Q\) of
		distinct patterns that occur satisfies \(Q\le 2^{4s^3}\).
		
			At each quadruple of positions, at least a \(p\)-fraction of all
		choices give a tetrahedron. Thus the average of \(|J(X)|/s^4\) is at
		least \(p\). As in~\eqref{eq:tri-exact}, the exact pattern count and
		AM--GM now give
		\[
		2^{-4s^3}\le Q^{-1}
		\le \left(\prod_{X\in\cX}\rho(J(X))\right)^{1/N^{4s}}
		\le \rho([s]^4)^p,
		\qquad\text{so}\qquad
		\rho([s]^4)\ge 2^{-4s^3/p}.
		\]
			Take \(s=2t\) and consider a full choice, so that every transversal
		of the four lists is a tetrahedron. If each list contains at least
		\(t\) distinct vertices, select \(t\) of them from each list.
		Every transversal of the selected sets is still a tetrahedron,
		giving a copy of \(K_4^{(3)}[t]\).
			We discard any full choice in which some list contains at most \(t-1\)
		distinct vertices. All \(2t\) entries of such a list lie in some
		\((t-1)\)-set in its class. 
		Fix a class \(V_i\) and a set \(A\subseteq V_i\) of size \(t-1\).
		The fraction of all choices for which every entry of the \(i\)th list
		lies in \(A\) is \(((t-1)/N)^{2t}\), since the other three lists
		are unrestricted. Every discarded full choice is counted for at least
		one pair \((i,A)\). Summing over these pairs, we conclude that the
		fraction of all choices that are full and survive is at least
		\[
		\rho([2t]^4)
		-4\binom{N}{t-1}\left(\frac{t-1}{N}\right)^{2t}.
		\]
		If no \(K_4^{(3)}[t]\) exists, the fraction of surviving choices is zero, so the preceding lower bound is at most zero. Hence, using
		the lower bound for \(\rho([2t]^4)\) and the estimate
		\(\binom{N}{t-1}\le (eN/(t-1))^{t-1}\), we obtain
		\[
		2^{-32t^3/p}
		\le 4\binom{N}{t-1}\left(\frac{t-1}{N}\right)^{2t}
		\le 4e^{t-1}\left(\frac{t-1}{N}\right)^{t+1}
		\le 4\left(\frac{e(t-1)}{N}\right)^t.
		\]
		Now choose \(t=\lfloor\sqrt{(p/64)\log_2 N}\rfloor\).
		Taking \(t\)th roots and multiplying by \(N\) gives
		\[
		\sqrt{N}\le N\,2^{-32t^2/p}\le 4^{1/t}e(t-1).
		\]
		The right-hand side is \(O(\sqrt{\log N})\), a contradiction for
		large \(N\).
		\end{proof}


\begin{thebibliography}{99}
		\small
		\setlength{\itemsep}{2pt}
		\setlength{\parskip}{0pt}
		
		\bibitem{ADLRY}
		N. Alon, R. A. Duke, H. Lefmann, V. R\"odl and R. Yuster,
		\href{https://doi.org/10.1006/jagm.1994.1005}{The algorithmic aspects of the regularity lemma},
		\emph{J. Algorithms} \textbf{16} (1994), 80--109.
		
		\bibitem{BE}
		B. Bollob\'as and P. Erd\H{o}s,
		\href{https://doi.org/10.1112/blms/5.3.317}{On the structure of edge graphs},
		\emph{Bull. London Math. Soc.} \textbf{5} (1973), 317--321.
		
		\bibitem{BES}
		B. Bollob\'as, P. Erd\H{o}s and M. Simonovits,
		\href{https://doi.org/10.1112/jlms/s2-12.2.219}{On the structure of edge graphs II},
		\emph{J. London Math. Soc.} (2) \textbf{12} (1976), 219--224.
		
		\bibitem{CFS10}
		D. Conlon, J. Fox and B. Sudakov,
		\href{https://doi.org/10.1090/S0894-0347-09-00645-6}{Hypergraph Ramsey numbers},
		\emph{J. Amer. Math. Soc.} \textbf{23} (2010), 247--266.
		
		\bibitem{CFS}
		D. Conlon, J. Fox and B. Sudakov,
		\href{https://doi.org/10.1007/s11856-011-0016-6}{Large almost monochromatic subsets in hypergraphs},
		\emph{Israel J. Math.} \textbf{181} (2011), 423--432.
		
		\bibitem{Erdos64}
		P. Erd\H{o}s,
		\href{https://doi.org/10.1007/BF02759942}{On extremal problems of graphs and generalized graphs},
		\emph{Israel J. Math.} \textbf{2} (1964), 183--190.
		
		\bibitem{Erdos90}
		P. Erd\H{o}s,
		\href{https://doi.org/10.1007/978-3-642-72905-8_2}{Problems and results on graphs and hypergraphs: similarities and differences},
		in \emph{Mathematics of Ramsey Theory}
		(J. Ne\v{s}et\v{r}il and V. R\"odl, eds.), Algorithms Combin. \textbf{5},
		Springer, Berlin, 1990, 12--28.
		
		\bibitem{EK}
		P. Erd\H{o}s and D. J. Kleitman,
		\href{https://www.renyi.hu/~p_erdos/1968-05.pdf}{On coloring graphs to maximize the proportion of multicolored \(k\)-edges},
		\emph{J. Combin. Theory} \textbf{5} (1968), 164--169.
		
		\bibitem{ES66}
		P. Erd\H{o}s and M. Simonovits,
		\href{https://www.renyi.hu/~p_erdos/1966-01.pdf}{A limit theorem in graph theory},
		\emph{Studia Sci. Math. Hungar.} \textbf{1} (1966), 51--57.
		
		\bibitem{ES83}
		P. Erd\H{o}s and M. Simonovits,
		\href{https://doi.org/10.1007/BF02579292}{Supersaturated graphs and hypergraphs},
		\emph{Combinatorica} \textbf{3} (1983), 181--192.
		
		\bibitem{ES46}
		P. Erd\H{o}s and A. H. Stone,
		\href{https://doi.org/10.1090/S0002-9904-1946-08715-7}{On the structure of linear graphs},
		\emph{Bull. Amer. Math. Soc.} \textbf{52} (1946), 1087--1091.
		
		\bibitem{FS}
		J. Fox and B. Sudakov,
		\href{https://doi.org/10.1002/rsa.20344}{Dependent random choice},
		\emph{Random Structures Algorithms} \textbf{38} (2011), 68--99.
		
		\bibitem{FWZ}
		J. Fox, Y. Wigderson and Y. Zhou, Finding blowups one vertex at a time,
		\href{https://arxiv.org/abs/2605.23301}{arXiv:2605.23301}, 2026.
		
		\bibitem{FK}
		H. Furstenberg and Y. Katznelson,
		\href{https://doi.org/10.1007/BF02790016}{An ergodic Szemer\'edi theorem for commuting transformations},
		\emph{J. Analyse Math.} \textbf{34} (1978), 275--291.
		
		\bibitem{KST}
		T. K\H{o}v\'ari, V. T. S\'os and P. Tur\'an,
		\href{https://doi.org/10.4064/cm-3-1-50-57}{On a problem of K. Zarankiewicz},
		\emph{Colloq. Math.} \textbf{3} (1954), 50--57.
		
		\bibitem{NikClique}
		V. Nikiforov,
		\href{https://doi.org/10.1112/blms/bdm093}{Graphs with many \(r\)-cliques have large complete \(r\)-partite subgraphs},
		\emph{Bull. Lond. Math. Soc.} \textbf{40} (2008), 23--25.
		
		\bibitem{NikGeneral}
		V. Nikiforov,
		\href{https://doi.org/10.37236/881}{Graphs with many copies of a given subgraph},
		\emph{Electron. J. Combin.} \textbf{15} (2008), Note 6, 6 pp.
		
		\bibitem{RS}
		V. R\"odl and M. Schacht,
		\href{https://doi.org/10.1002/rsa.20441}{Complete partite subgraphs in dense hypergraphs},
		\emph{Random Structures Algorithms} \textbf{41} (2012), 557--573.
		
	\end{thebibliography}
\end{document}